\documentclass[10pt,a4paper]{article}
\usepackage{amsmath}\usepackage{theorem}\usepackage[all]{xy}
\usepackage{amssymb}\usepackage{latexsym}

\newtheorem{The}{Theorem}
\newtheorem{Lem}[The]{Lemma}\newtheorem{Pro}[The]{Proposition}

\newcommand{\bevis}{\textit{Proof: }}
\newcommand{\qed}{\hfill$\Box$}
\newcommand{\ner}{\vspace{\theorempostskipamount}}
\newcommand{\ra}{\rightarrow}\newcommand{\hra}{\hookrightarrow}

\newcommand{\CP}{\mathit{CP}}
\newcommand{\R}{\mathbf{R}}
\newcommand{\C}{\mathbf{C}}
\newcommand{\Z}{\mathbf{Z}}
\newcommand{\pnt}{{\{\mathit{pt.}\}}}
\newcommand{\Top}{\mathrm{Top}}

\newcommand{\stj}{\widetilde{*}}

\begin{document}
\title{The linear isometries operad \\ in Lie--Tate homology}

\author{\textsc{Pelle Salomonsson}}\date{\small{\texttt{pellegestalten@gmail.com}}}
\vspace{.1cm}
\maketitle
\vspace{.5cm}
\abstract{}
\!\!\!\!\!\!\!\noindent
\begin{quote} We give an independent, and perhaps somewhat simplified, description of the product in negative Tate-cohomology (the generalised version for compact Lie-groups). We describe, but do not compute, the corresponding action of the Dyer--Lashof-algebra, using the linear-isometries operad.
\end{quote}

\vspace{.7cm}
\noindent We give a description of the ``join-product,'' using space-level constructions (as opposed to chain-complexes). Then we describe the corresponding operad-action, using the linear-isometries operad.

I came up with this construction independently---and only later realised that it seems to be (in some sense) implicit in the general machinery~\cite{mcclure} of equivariant stable homotopy-theory.

\section{The join-product}
We fix a space $X$ and consider the category $\Top_X$ of spaces mapping to it. This category can be given a symmetric monoidal structure, as follows.
Consider the space $\Delta^1\!\times\!\mathrm{Map}(\Delta^1\!,\,X)$ of ``pointed singular one-simplices'' in~$X$. Let $\mathrm{Simp}_1X$ be the quotient-space given by making edge-identifications: that is, if the distinguished point lies at one of the vertices, then we forget everything about the mapping except the image of that point.

Next let two spaces $Y_1,Y_2$ mapping to~$X$ be given. Given a pointed singular one-simplex, assume given also (for $i=1,2$) a choice of a point in~$Y_i$ mapping to the image of the $i$'th vertex of $\Delta^1$. Then if the distinguished point lies at one the vertices: forget about everything except the point in the corresponding $Y_i$. Let the space of all such mappings and all such choices be denoted by $Y_1*_XY_2$, or sometimes simply by $Y_1*Y_2$. In fact, \emph{throughout in this exposition},~$*$ denotes~$*_X$. Evaluation at the distinguished point gives a mapping to~$X$. We separately legislate that $\varnothing*Y=Y*\varnothing=Y$.

\begin{Pro} The product $*_X$ makes the homotopy-category of $\Top_X$ into a monoidal category, having~$\varnothing$ as identity-object.
\end{Pro}
\bevis Let $\mathrm{Simp_2}X$ be defined analogously to $\mathrm{Simp_1}X$. There is a rather evident mapping $\mathrm{Simp_2}X\ra X\!*\!(X\!*\!X)$, which in fact can be seen to be a homotopy-equivalence (by a homotopy that does not move the vertices of~$\Delta^2$). The associativity-claim follows from this fact, and from the independence of $\mathrm{Simp_2}X$ with respect to the choice of bracket-pattern. \qed

\ner\noindent Let a point in $Y*Z$ be given by $f:\pnt\ra\Delta^1$ and $g:\Delta^1\ra X$, where of course~$f$ is not a simplicial map. Then we have a
natural involution given by postcomposing~$f$ and precomposing~$g$ with the flip $\Delta^1\ra\Delta^1$. It makes (the homotopy-category of) $\Top_X$ into a symmetric monoidal category.

\subsection{The product in homology}
Let us restrict attention to the case when~$X$ is the classifying space of a compact Lie-group, when a refined model of the join-product becomes available. Namely, the free pathspace-fibration $\Omega BG\hra PBG\ra BG^{\times2}$ admits a deformation-retraction onto a subbundle with fibre homeomorphic to~$G$.

Let~$H_*$ denote ordinary homology with mod~$p$ coefficients, and disregard the (completely uninteresting) case when~$G$ is finite group of order prime to~$p$. As the space  $\Sigma\Omega BG\cong\Sigma G$ then has a welldefined fundamental-class, remaining well-defined in families, one readily constructs the following maps, essentially of integration-along-fibres type:

\begin{Pro}\label{umkehrproposition} Assume that $G$ is not a finite group of order prime to~$p$. For all~$n$ there are then natural maps $H_*(BG^{\times n})\ra  H_*(BG^{*n})$.
\end{Pro}
\bevis We need to make use of the non-associative product~$\stj$ given by forgetting slightly less than we did before. When the distinguished point lies at one of the vertices: remember \emph{both} of the chosen points in the $Y_i$'s. There are forgetful mappings $(Y_1\stj\cdots)\ra (Y_1*\cdots)$ for any bracket-pattern. The advantage of~$\stj$ is that there is a fibrebundle
$(Y_1\stj\cdots)\ra X^{\times n}$ with fibre a certain twisted direct product of copies of $\Sigma\Omega X$ (a tower of fibrations). It remains to prove that this construction is independent of the chosen bracket-pattern. But any two different umkehr-maps up to two different fibrations mapping to some common third space by fibrewise compatible ``generic isomorphisms'' must in fact (at least up to sign) agree when pushed-forward to the homology of that space. The assumption on~$G$ is need to get a closed chain, as we are forced to work with the unreduced suspension.
\qed

\ner\noindent For any $X$, regarded as an object of $\Top_X$, the structure-map $X^{*n}\ra X$ is in fact a homotopy-equivalence. We get then for each~$n$, using the K\"{u}nneth-isomorphism, the following composed map.
\begin{equation}\label{homologymultiplication}
\xymatrix@R-0pc@C-0pc{
   H_*(BG)^{\otimes n}  \ar[r] & H_*(BG^{*n}) \ar[r]^{=} & H_*(BG)
   }\,\,.
\end{equation}
It increases homological degree by $(n\!-\!1)(\dim\,G\!+\!1)$. It should not be confused with the Pontrjagin-multiplication on loopspaces---which is also, by a pure coincidence of notation, denoted by~$*$. For some further discussion (based on the work of Kreck), and for the connection to Tate-cohomology, see the manuscript~\cite{tene} of H.~Tene.

\begin{Pro} The join-product is graded-commutative in the sense that $a*b=(-1)^{s}b*a$, where $s=\deg\,a\cdot\deg\,b+\dim\,G+1$. In particular, $\pnt^{*2}=[\Sigma G]=0$ if~$p$ is odd and $\dim\,G$ is even.
\end{Pro}
\bevis One easily reduces to the case when $\deg a=\deg b=0$. Hence now let~$a$ and~$b$ denote two distinct points in~$X$. The space $P_{(a,b)}X$ of paths going from~$a$ to~$b$ has a certain orientation (remaining globally  welldefined as~$a$ and~$b$ vary). It gives an orientation to the space $P_{\{a,b\}}X$ of unoriented paths going between these points. Indeed, if we (non-canonically) identify  $P_{(a,b)}X$ with~$G$, then we get an identification $P_{\{a,b\}}X=G$ as well. However, we get another such identification by using $P_{(b,a)}X$ instead, and the two identifications can be seen to be related through the inverse map $G\ra G :g\mapsto g^{-1}$, which is orientation-reversing when $\dim G$ is odd. This essentially proves the claim, upon thinking through how the integration-along-fibres procedure is defined. Of course, the extra ``$\,+1\,$'' appearing in the formula comes from the flip of~$\Delta^1$.   \qed

\section{The linear-isometries operad}
We prove first a lemma. For $W$ a constant $X$-space (that is, mapping to a single point in~$X$) we put $W\times Y$ to be $W\times Y\ra Y\ra X$. We have then:
\begin{Lem} \label{movinglemma} Let $Y$, $Z$ and $W$ be $X$-spaces, with $W$ constant. Then there is a natural map $W\times(Y*Z)\ra (W\times Y)*Z$.
\end{Lem}
\bevis The map forgets the $W$-coordinate whenever the distinguished point lies at the $Z$-endpoint.   \qed

\ner\noindent Let $L_*$ be the so-called linear-isometries-operad, so that $L_n$ para\-meterises linear embeddings $\bigl(\R^{\infty}\bigr)^{\times n}\hra\R^{\infty}$. It is a $\Sigma$-free contractible operad. Let an inclusion
\begin{equation}\label{linearrepresentation}
\xymatrix@R+0pc@C-.7pc{
   G\,\, \ar@{^(->}[r] & \,\,O(k)
   } \vspace{.2cm}
\end{equation}
of~$G$ into some finite-dimensional rotation-group be given. Then consider the corresponding Stiefel-space of ordered orthogonal vectors in~$\R^\infty$\!. It is a contractible space, upon which~$G$ acts freely. The quotient is~$BG$.

The space $BG$ admits embeddings of homeomorphic copies of itself into itself. Indeed, such an embedding can be obtained from an embedding $\R^\infty\hra\R^\infty$. Given a set of such subspaces, disjoint from each other,  we may form the join of them. The union of these joins make up a fibrebundle over~$L_n$ with fibre $BG^{*n}$ (regarded now as an abstract space, without canonical mapping to $BG$). By associating points that map identically to  $BG$ we equip the bundle with an equivariant structure, acting nontrivially on insertions. Moreover, as the base is contractible the bundle can be (noncanonically) trivialised. We get a $\Sigma_n$-equivariant map
\begin{equation}\label{bundletrivialisation}
\xymatrix@R+0pc@C-0pc{
   L_n\times BG^{*n}  \ar[r] & BG
   }\,\,, \vspace{.2cm}
\end{equation}
which, if we regard~$L_n$ as a constant $BG$-space, can be considered a map in $\Top_{BG}$.

Various notions of ``bimonoidal categories'' have been studied. We would say that a category is a ``priority-swapping category'' if it is equip\-ped with a distinguished subcategory and two symmetric monoidal structures satisfying the conclusion of lemma~\ref{movinglemma}, together with some further axioms involving associativity-isomorphisms. In such a category we would concern ourselves with operads with respect to the product of $\times$-type, and study bimonoidal actions making diagrams of the following type commute equivariantly (where we for ease of exposition confine ourselves to the $\Sigma_2\wr\Sigma_2$-case):
\begin{equation}\label{operaddiagram}
\xymatrix@R-1.4pc@C-.9pc{
   L_2\times \Bigl((L_2\times X^{*2})*
   (L_2\times X^{*2})\Bigr) \ar[r]
   &L_2\times (X*X)\ar[dr]&\\
   &&X\\
   L_2\times\Bigl(L_2\times L_2\times \bigl(X^{*2}*X^{*2}\bigr)\Bigr)
   \ar[uu]_{\text{Lemma~\ref{movinglemma}}}\ar[r]
   &L_4\times X^{*4}\ar[ur]&
   }
\end{equation}
\begin{Pro} With $L_*$the linear-isometries operad and $X=BG$, diagram~\ref{operaddiagram} and its analogues commute equivariantly on the nose.
\end{Pro}
\bevis True, the splitting~(\ref{bundletrivialisation}) is canonical only up to local slidings along fibres by the $G^{}$-action. However, the spaces and the mappings appearing in the diagram are ``independent'' of the splitting: another choice gives an isomorphic diagram. The statement follows from the naturality of the definitions. \qed

\ner\noindent We write $c_i$ for some nonclosed chain in the chain-complex of $L_p$ mapping to the standard-generator of $H_i\bigl((L_p)_{\Z/p}\bigr)$. As for the meaning of ``standard,'' we think of the unique cell in the standard-model of $B(\Z/p)$. If $p\!=\!2$, there is of course no ambiguity. Next, given a homology-class~$y$ in $BG$, we get a homology-class in $BG^{\times p}$, and then by proposition~\ref{umkehrproposition} a~homology-class $y^{*p}$ in $BG^{*p}$. We proceed to consider the class
$c_i\otimes_{\Sigma_p} y^{*p}$ in the equivariant homology of $L_p\times BG^{*p}$. We write $\theta_i(y)$ for its image under~(\ref{homologymultiplication}) in $H_{*}BG$. In the classical case it transpires~\cite{dyer-lashof} that most of the $\theta_i$'s vanish when~$p$ is odd. This holds true in our setting as well, for the same reason.
Thus, when~$p$ is odd we write $Q_i:=\theta_{2i(p-1)}$.
In the $p\!=\!2$ case we write $Q_i:=\theta_i$.

\begin{The}\label{dyerlashofaction} The action of the operators $Q_i$, acting on $H_*BG$, is in fact an action of the (lower-index) Dyer--Lashof-algebra---that is, the Adem-relations (as well as the Cartan-formula with respect to the $*$-product) lie in its kernel.
\end{The}
\bevis We must carefully think through what diagram~\ref{operaddiagram} means in terms of iterated operations, as the unusual monoidal structure might conceivably be a hindrance to such an interpretation. Let $c_i$, $c_j$ and $y$ be given, with notation as above. In the upper righthand space, let us consider the class $c_i\otimes (Q_jy)^{*2}$. It is, by the naturality of the $z\mapsto z^{*2}$ transformation, the image of the class $c_j\otimes (c_i\otimes y^{*2})^{*2}$\!, which in its turn is the image of $c_j\otimes c_i^{\otimes2}\otimes y^{*4}$. Hence beginning with this class in the lower lefthand corner and following the upper path, we indeed get the iterated action on~$y$. On the other hand the lower path gives us, more or less by definition, the Adem-relations. The same reasoning applies in the case of odd~$p$, as well as in the case of the Cartan-formula.
\qed

\ner\noindent \textbf{Remark I.} Tate-cohomology has been generalised in two directions: to the case when~$G$ is a compact Lie-group, and to general $G$-spaces (``Swan cohomology''). In this paper we say nothing about the spatial version.

\ner\noindent \textbf{Remark II.} We also say nothing about the question as to whether the action is independent of the choice of linear representation $G\hra O(k)$.

\ner\noindent \textbf{Remark III.} In the case of a finite~$G$, an action of the Dyer--Lashof-algebra has been constructed~\cite{langer} without recourse to any linear representation. There is also a construction by McClure~\cite{mcclure}. That operad acts, however, on the spectrum that represents the cohomology-theory, not on the space whose cohomology we wish to compute.

\ner\noindent \textbf{Remark IV.} In the case of a higher-dimensional connected~$G$, the product seems to vanish identically (except in some low-dimensional cases). This is claimed in a remark in~\cite{greenlees-may}. One would perhaps expect the same to hold true for the operations. And in that case, the developments in this paper seem to be of interest mostly for finite groups.

\section{The case $G=S^1$}
The action on $H_*BG$ for $G=\Z/p$ has been computed by M.~Langer~\cite{langer}. It seems clear that his action, at least in this case, coincides with ours. Hence, as for the abelian case, it essentially remains for us to deal with the case $G=S^1$\!. Write~$x_i$ for the standard-generator of $H_{2i}(\CP^\infty\!\!,\,\Z/p)$. We state the result only in the $p=2$ case: these direct cellular techniques seem unworthwhile in general.

\begin{The} Assume $p=2$. Then $Q_{2j}(x_i)=\binom{i+j}{j}x_{2i+j+1}$.
\end{The}
\bevis One proves that  $Q_{2j}(x_0)=x_{j+1}$, and then the rest follows from the Cartan-formula. To prove this, we must go back to the linear-isometries operad and construct a family of embeddings $(\R^\infty)^{\times2}\ra\R^\infty$, parameterised over a $2j$-dimensional cell~$c$ which becomes a nontrivial cycle in the quotient $(L_2)_{\Z/2}$. We endow $\R^\infty$ with its standard complex structure. Moreover, it suffices to construct the family on $(\C^1)^{\times2}\subset(\C^\infty)^{\times2}$. Well then, consider the standard inclusion $S^{2j}\hra\C^{j+1}$, and identify $\C^{j+1}$ with any affine hyperplane $V\subset\C^{j+2}$ not passing through the origin. Sending the first basis-vector of $(\C^1)^{\times2}$ to any upper-hemisphere point and the second to its antipode gives the desired family of linear embeddings. We get then, under the quotient by the $\C^*$-action, a family of lines $\CP^1=\pnt^{*2}\subset\CP^{j+1}$ parameterised over $\mathit{RP}^{2j}$. One verifies that there passes exactly one line through a generic point on $\CP^{j+1}=x_{j+1}$, which was we wanted to prove. \qed


\end{document}